\documentclass[11pt,a4paper,twoside]{article}
\usepackage{latexsym,amsfonts,amsmath, amsthm, amssymb}
\usepackage{fullpage}
\usepackage{xcolor,bm, bbm}
\usepackage[utf8x]{inputenc}
\usepackage{array}
\usepackage{mathrsfs}
\usepackage{tikz}
\usepackage{pgfplots}

\usepackage{fourier}

\newcommand{\R}{\mathbb R}
\newcommand{\N}{\mathbb N}

\newcommand{\E}{\mathbb E}
\newcommand{\Pro}{\mathbb P}

\def\dint{\textup{d}}

\newtheorem{thmalpha}{Theorem}

{
\theoremstyle{definition}

}

\allowdisplaybreaks
\begin{document}

\title{\bf Talagrand's mathematical journey to the Abel Prize 2024} 

\medskip

\author{Olivier Gu\'edon and Joscha Prochno}



\date{}

\maketitle

\begin{abstract}
\small Michel Talagrand (Centre National de la Recherche Scientifique, France) has been awarded the prestigious Abel Prize for 2024
for his work in probability theory, functional analysis, and statistical physics. In this note, we introduce the Abel Prize Laureate and his main contributions.

\medspace
\end{abstract}
 
\tableofcontents

\section{Introduction}

\emph{``If I had heard that an alien mothership had landed in front of the city hall, I don't think I would've been more surprised. I had never, ever thought that this could happen.''} was the humble reaction of Michel Talagrand from the French National Centre for Scientific Research, (CNRS) when he was awarded the prestigious 2024 Abel Prize coming with a monetary award of 7.5 million Norwegian kroner (about $\$$700,000). The Norwegian Academy of Science and Letters had decided that Talagrand should receive the prize ``for his groundbreaking contributions to probability theory and functional analysis, with outstanding applications in mathematical physics and statistics.'' and the authors of this note could not be happier about the decision to award this extraordinary and outstanding mathematician who has published 299 papers as well as several books (by May 2024), and has had (and still has) a strong influence on our discipline. Three areas of mathematics, and major playing fields of Talagrand, are highlighted in the announcement, namely (A) suprema of stochastic processes, (B) concentration of measures, and (C) spin glasses. After a short biography, we shall provide a glimpse at Talagrand's contributions to those three areas.

\subsection{Biographical Sketch}

Michel Talagrand was born 15 February 1952 in B\'eziers and grew up in Lyon (France), where he later studied mathematics. Already his father Pierre was a high school teacher of mathematics and his mother Raymonde a French teacher in junior high school. His interest in the sciences was sparked by the scientific magazine ``Sciences et Avenir'' to which his father had subscribed. In his short biography for the Shaw Prize, Talagrand describes himself as a ``very average student'', how his father tried to get him interested in mathematics (teaching him basic group theory) and how, after already having lost sight in his right eye due to retinal detachment aged only 5, getting a retinal detachment in the left eye at age 14 (even though it got treated) changed the course of his life and made him a different person. He had to spend a month in the hospital and his father came each day to teach Talagrand, his eyes bandaged, mathematics; this is how he learned the power of abstraction. In his last year of high school, he ranked third in both math and physics in the French Olympiad for the top high school seniors. Talagrand continued to study mathematics and physics at Universit\'e de Lyon and, after two years, chose mathematics as his major subject. After his master's degree he took part in the very high level french competition named ``Agrégation de Mathématiques'' and came out first with exceptional grades. The very same year, he applied at CNRS and was accepted in 1974. Talagrand then moved to Paris where he continued his life-long course of being a mathematician, which, among others, has been strongly influenced by Gilles Pisier (born 18 November 1950) and Vitali D. Milman (born 23 August 1939). For more personal information about Michel Talagrand, the reader may have a look at the interview (in french) that Talagrand gave to Gilles Godefroy (a mathematician and close friend) for the ``Gazette des Mathématiciens'' \cite{MR3965513}.

In 1977, Talagrand completed his PhD at Pierre and Marie Curie University (Paris VI) with the thesis  \emph{Mesures invariantes, compacts de fonctions mesurables et topologie faible des espaces de Banach} written under supervision of Gustave Choquet (1 March 1915 -- 14 November 2006). From 1985 until his retirement in 2017, he had been Directeur de Recherches at CNRS and a member of the "\'Equipe d'Analyse" which became the Functional Analysis Team of the Institut de Mathématiques de Jussieu (IMJ-PRG). Talagrand was elected as fellow of the French Scientific Academy in 2004 and was invited speaker at the International Congress of Mathematicians in Kyoto in 1990 (Some isoperimetric inequalities and their application) and in Berlin in 1998 (Huge Random Structures and Mean Field Model of Spin Glasses). While for most mathematicians being awarded \emph{only} one of the major prizes in mathematics would already be considered a knightly accolade, Talagrand has received various of them, for instance, the Lo\`eve Prize (1995), the Fermat prize (1997), Shaw Prize (2019), Stefan Banach Medal (2022), and now the Abel Prize (2024).

\section{Talagrand's Work -- Probability in Banach spaces}
 
 Let us now continue with short discussions of some of Talagrand's contributions and works winning him the Abel Prize. 
 
 \subsection{Suprema of stochastic processes}
 
 The fundamental question that underlies one of the major contributions of Talgrand is one that can be rather easily formulated and understood. Say we are given a collection of random variables $(X_t)_{t\in T}$ for some index set $T$. The goal is to understand the supremum of such a stochastic process, more precisely, to find upper and lower bounds for the quantity 
   \begin{equation}\label{eq:expected sup definition}
     \E \sup_{t\in T} X_t :=  \sup \Big\{ \E \sup_{t\in F} X_t \,:\, F\subseteq T, \, F \text{ finite} \Big \}
   \end{equation}
i.e., for the expected supremum of the process; the definition above takes into account that the set $T$ may be uncountable. This also shows that the crucial case for estimating the expected supremum is the one of finite $T$. 

The arguably most natural and important case is the one of (mean-zero) Gaussian processes and studying under which conditions those processes are bounded and understanding their asymptotic behavior has a long tradition. Let us recall that a random process $(X_t)_{t\in T}$ is called Gaussian process if and only if, for any finite $F\subseteq T$, the random vector $(X_t)_{t\in F}$ has a normal distribution or, equivalently, if every finite linear combination $\sum_{t\in F} a_tX_t$, $a_t\in\R$, is a normal random variable; we shall only be interested in Gaussian processes of mean zero here. It is not hard to show that the increments 
  \begin{equation}\label{eq:canonical metric}
    d(s,t):= \| X_s - X_t \|_2 = \big(\E[\big|X_s-X_t\big|^2]\big)^{1/2}, \quad s,t\in T,
  \end{equation}
 of the random process always define a (canonical) metric on the index set $T$, even if this is an abstract set having no geometric structure, and that they also determine the covariance function and thus the distribution of the process. This means that questions about the distribution of the process $(X_t)_{t\in T}$ are intimately connected to the geometry of the metric space $(T,d)$. So in order to study upper and lower bounds for the overall magnitude of the process, i.e., $\E \sup_{t\in T} X_t$, we may study how this magnitude relates to the geometry of $(T,d)$? 
 
 As it is quite instructive and provides clear motivation, we shall start with two classical results of Sudakov and Dudley. Even though the bounds are not optimal, one can see how the processes extremal properties are related to geometric characteristic of the metric space. More precisely, we shall see an intimate connection between the expected supremum of the Gaussian process and the notion of metric entropy. Let us recall the essential definitions here. Given a metric space $(T,d)$, a subset $A\subseteq T$, and $\varepsilon\in(0,\infty)$, we call $\mathcal N\subseteq A$ an $\varepsilon$-net of $A$ if and only if 
   \[
     \forall x\in A \,\, \exists x_0\in\mathcal N:\quad d(x,x_0) \leq \varepsilon.
   \]
We call the smallest possible cardinality of an $\varepsilon$-net of $A$ the \emph{$\varepsilon$-covering number} of $A$ and write $N(A,d,\varepsilon)$ (note that this number may be $\infty$ if $T$ does not admit a finite $\varepsilon$-net). Equivalently, $N(A,d,\varepsilon)$ is the smallest number of closed metric-balls with centers in $A$ and radii $\varepsilon$ whose union covers the set $A$; for those more familiar with packing numbers, we note that those concepts are essentially equivalent. We can now define the \emph{metric entropy} of $(T,d)$ as the logarithm of the covering number, i.e., as $\log_2\big(N(T,d,\varepsilon)\big)$. 

The following result (see, e.g., \cite[Theorem 3.18]{LedouxTalagrand1991}) is commonly referred to as Sudakov minoration inequality, named after the first such result stated in \cite[Theorem 2]{Sudakov1971}. 
 
\begin{thmalpha}[Sudakov minoration inequality]
Let $(X_t)_{t\in T}$ be a mean-zero Gaussian process. Then, for any $\varepsilon\in(0,\infty)$, we have
  \[
    \E \sup_{t\in T} X_t \geq c \varepsilon \sqrt{\log_2\big(N(T,d,\varepsilon)\big)}
  \]
for some absolute constant $c\in(0,\infty)$, where $d$ is the canonical metric associated with the process as defined in \eqref{eq:canonical metric}. 
\end{thmalpha}

As is shown in detail in \cite{MR298726, Fernique1975, Sudakov1976}, the result can be obtained from a famous comparison result for Gaussian processes due to Slepian \cite[Lemma 1]{Slepian1962}, a result which was later generalized by Fernique \cite{Fernique1975}. This
shows that if for two mean-zero Gaussian processes one dominates the other in terms of the expectation of the second moment of increments, then this dominance carries over to the expected supremum of the processes. Note that if the metric space is not compact, then the covering numbers are infinite and so is the supremum of the process.

Let us continue with an upper bound that is also given in terms of covering numbers and referred to as Dudley's integral inequality as it can be traced back to ideas of Dudley from 1967 \cite{Dudley1967}; indeed, he performed all essential mathematical steps, but did not state the result explicitly. We shall only state the result for Gaussian processes, even though we could allow, e.g., processes with sub-Gaussian increments. We refer the reader to \cite[Proposition 2.2.10]{Talagrand2014} or \cite[Theorem 8.1.3]{Vershynin2018}.

\begin{thmalpha}[Dudley's integral inequality]
Let $(X_t)_{t\in T}$ be a mean-zero Gaussian process. Then  
\[
    \E \sup_{t\in T} X_t \leq C \int_0^{\infty}  \sqrt{\log_2\big(N(T,d,\varepsilon)\big)} \,\dint \varepsilon
  \]
for some absolute constant $C\in(0,\infty)$, where $d$ is the canonical metric associated with the process as defined in \eqref{eq:canonical metric}. 
\end{thmalpha}

The previous result can be obtained from a rather basic version of the chaining method, a powerful tool that is a sort of multi-scale version of the classical $\varepsilon$-net argument, where one builds progressively finer approximations to the points with finer nets. Dudley's bound is not always sharp and in fact those cases where it is not cannot be dismissed as exotic shown, for instance, by the case of ellipsoids of a Hilbert space (see \cite[Section 2.5]{Talagrand2014}). Let us also remark that in the two bounds we have seen so far, i.e., Sudakov's minoration inequality and Dudley's integral inequality, there is a clear gap, because if we look at $\sqrt{\log_2\big(N(T,d,\varepsilon)\big)}$ as a function in $\varepsilon$, then essentially, while Dudley's inequality bounds the expected supremum of the process by the area under the curve, Sudakov's inequality bounds it from below by the largest area of a rectangle under the curve. The reason behind all this is that the covering numbers do not contain enough information to control the magnitude of the process; an example showing this is given in the groundbreaking paper \cite[p.~101--102]{Talagrand1987} of Talagrand. Fernique apparently conjectured as early as 1974, after obtaining some improvement upon Dudley's result, that the existence of so-called majorizing measures (probability measures representing the most general possible form of a chaining argument going back to Kolmogorov) might characterize the boundedness of Gaussian processes instead. However, since the concept of majorizing measures was considered exotic by many others and had the reputation of being opaque and mysterious, Fernique remained rather isolated in his attempt to obtain a characterization. 

The continuation of the story, and behind his crowning achievement \cite{Talagrand1987}, is that in 1983 Gilles Pisier encouraged him to work on Fernique's conjecture and in fact kept goading him until he proved his famous \emph{majorizing measure theorem} in 1985 in which he proved the existence of a majorizing measure in a non-constructive way. Let us note that the proof of this result is considered notoriously difficult. In 1992 \cite{Talagrand1992}, Talagrand gave a completely elementary proof of the existence of majorizing measures for bounded Gaussian processes relying upon Sudakov minoration, the concentration of measure phenomenon, and a rather simple construction. Later, around 2000, Talagrand himself eschewed the use of majorizing measures in favor of a purely combinatorial approach referred to as \emph{the generic chaining} (see, e.g., \cite{Talagrand2014}).

In order to present Talagrand's majorizing measure theorem, which we formulate in its equivalent $\gamma_2$-functional or generic chaining version, we require some preparation; we refer to \cite{Talagrand2014} for a discussion of majorizing measures versus generic chaining, in particular to \cite[Section 6.2]{Talagrand2014} to explain the name majorizing measure theorem. Consider a metric space $(T,d)$. Then we call a sequence $(T_k)_{k=0}^\infty$ of subsets of $T$ an \emph{admissible sequence} (or just \emph{admissible}) if and only if the cardinalities satisfy
  \[
    |T_0|=1 \qquad\text{and}\qquad |T_k| \leq 2^{2^k} \text{ for $k\in\N$}.
  \]
The $\gamma_2$-functional of the metric space $(T,d)$ is then defined as 
  \[
    \gamma_2(T,d) := \inf_{(T_k)_{k=0}^{\infty}\text{ admissible}} \sup_{t\in T} \sum_{k=0}^\infty 2^{k/2} d(t,T_k).
  \]
It is not hard to check (see, e.g., \cite{Talagrand2014} or \cite[Chapter 8]{Vershynin2018}) that this expression is indeed smaller than Dudley's bound, which may be equivalently formulated in a similar fashion but with the supremum inside the sum on the right-hand side.

\begin{thmalpha}[Talagrand's majorizing measure theorem]
Let $(X_t)_{t\in T}$ be a mean-zero Gaussian process. Then  
  \[
    c\gamma_2(T,d)\leq \E \sup_{t\in T} X_t \leq C \gamma_2(T,d)
  \]
for some absolute constants $c,C\in(0,\infty)$, where $d$ is the canonical metric associated with the process as defined in \eqref{eq:canonical metric}. 
\end{thmalpha}

The upper bound in the previous theorem is an improvement upon Dudley's bound and commonly referred to as the \emph{generic chaining bound}. The major contribution of Talagrand was to prove the corresponding lower bound. Together the result shows us that chaining suffices to explain the size of a mean-zero Gaussian process and so we shall briefly sketch the chaining argument here. 

\begin{proof}[Proof of the generic chaining bound]
The definition of the expected supremum of a process in \eqref{eq:expected sup definition} allows us to assume w.~l.~o.~g.~that the index set $T$, which we endow with its canonical metric, is finite. We consider an admissible sequence $(T_k)_{k=0}^\infty$ of subsets of $T$ and let $T_0 := \{t_0 \}$; note that since the process is centered, we have 
  \[
    \E\sup_{t\in T} X_t = \E\sup_{t\in T} (X_t-X_{t_0})
  \]
and so we may work with the expected supremum of the displacement $X_t-X_{t_0}$, $t\in T$; in particular, the considered random variables $\sup_{t\in T} (X_t-X_{t_0})$ are non-negative now and so 
  \begin{equation}\label{eq:integral expression expectation}
    \E\sup_{t\in T} X_t = \E\sup_{t\in T} (X_t-X_{t_0}) = \int_{0}^\infty \Pro\big[\sup_{t\in T} (X_t-X_{t_0})>u\big]\,\dint u
  \end{equation}
is the expression we seek to bound from above.  

For a point $t\in T$, we shall denote by $\pi_k(t)\in T_k$ one of the closest points to $t$ in the set $T_k$ and by $\pi_k:T\to T_k$ the corresponding closest point mapping, i.e., 
  \[
  d(t,T_k) = d(t,\pi_k(t));
  \]
in other words, $\pi_k(t)$ is the best approximation to $t$ in $T_k$. Now we consider the chain taking us from $t_0$ to $t$ through a collection of admissible sets, more precisely, for $N\in\N$ large enough we walk along the chain
  \[
    t_0=\pi_0(t) \rightarrow \pi_1(t) \rightarrow \pi_2(t) \rightarrow \cdots \rightarrow \pi_N(t)=t;
  \]
in fact, since $T$ is finite, for some large enough $N_0\in\N$, we have $T_k=T$ for $k\geq N_0$. 
We are now able to decompose the increments along our chain. The idea behind this strategy is as follows: doing a union bound at the very beginning, i.e.,
  \[
    \Pro\big[\sup_{t\in T} (X_t-X_{t_0})>u\big] \leq \sum_{t\in T}\Pro\big[X_t-X_{t_0}>u\big],
  \]
is doomed to fail if we have many rather correlated random variables, but it should be rather effective along the chain of approximations in the canonical metric, because (A) there are fewer terms which are rather different and (B) the displaced random variables along the chain are smaller than the original displaced ones, making their supremum easier to control.

So let us express the displacement $X_t-X_{t_0}$ as a telescoping sum, i.e., 
  \[
    X_t-X_{t_0} = \sum_{k=1}^N \big(X_{ \pi_k(t)} - X_{\pi_{k-1}(t)}\big),
  \]
and try to control the increments in the sum, seeking to obtain a uniform bound 
  \[
    \big| X_{ \pi_k(t)} - X_{\pi_{k-1}(t)} \big| \leq 2^{k/2} d(t,T_k) \quad \forall k\geq 1  \,\, \forall t\in T
  \]
with high probability.

Let us recall the classical Gaussian concentration bound, telling us that for all $a,b\in T$ and any $\lambda>0$,
  \[
    \Pro\big[X_a-X_b > \lambda\big] \leq e^{-\frac{\lambda^2}{2d(a,b)^2}} ;
  \]
after all, $X_a-X_b$ is a mean-zero Gaussian with variance $d(a,b)^2$.  

Using this estimate in our setting, for all $k\geq 1$, $t\in T$, and $u>0$, we can deduce that 
  \begin{equation}\label{eq:gaussian concentration bound}
      \Pro\Big[ \big| X_{ \pi_k(t)} - X_{\pi_{k-1}(t)} \big| > u 2^{k/2} d(\pi_k(t),\pi_{k-1}(t)) \Big] \leq e^{-u^22^{k-1}}.
  \end{equation}
  
In order to unfix the parameter $t$ (for $k$ still being fixed, we shall later just use a standard union bound), we observe that the number of pairs $(\pi_k(t),\pi_{k-1}(t))$ is bounded by
  \[
    |T_k||T_{k-1}| \leq |T_k|^2 = 2^{2^{k+1}},
  \]
which allows us to obtain 
  \begin{equation}\label{eq:estimate unfixing t}
    \Pro\Big[ \exists t\in T\,:\, \big| X_{ \pi_k(t)} - X_{\pi_{k-1}(t)} \big| > u 2^{k/2} d(\pi_k(t),\pi_{k-1}(t)) \Big]  \leq 2^{2^{k+1}}e^{-u^22^{k-1}}.
  \end{equation}
For $u>0$, we now define the event
  \[
    \Omega_u := \Big\{ \forall k\geq 1\,\forall t\in T\,:\, \big| X_{ \pi_k(t)} - X_{\pi_{k-1}(t)} \big| \leq u 2^{k/2} d(\pi_k(t),\pi_{k-1}(t)) \Big\}.
  \]
Then, using a simple union bound in $k$ and summing the esimate in \eqref{eq:estimate unfixing t}, we obtain  
  \begin{equation}\label{eq:bounding omega complement}
    \Pro[\Omega_u^c] \leq \sum_{k\geq 1} 2^{2^{k+1}} e^{-u^2 2^{k-1}} \leq O(1)e^{-u^2}
  \end{equation}
for $u\geq c$, with $c\in(0,\infty)$ being some absolute constant ($c=4$ suffices). Let us define the quantity 
  \[
    S:= \sup_{t\in T} \sum_{k\geq 1} 2^{k/2}d(\pi_k(t),\pi_{k-1}(t)).
  \]
It is immediately clear that for any $\omega\in\Omega_u$, it holds that $\sup_{t\in T}(X_t-X_{t_0})(\omega) \leq u S$. Combining this with \eqref{eq:bounding omega complement}, we see that
  \[
    \Pro\big[ \sup_{t\in T}(X_t-X_{t_0}) > u S\big] \leq O(1) e^{-u^2}.
  \]
Using \eqref{eq:integral expression expectation} and keeping in mind that the integrand is less than or equal to $1$, this implies
  \[
    \E\sup_{t\in T}X_t = \E \sup_{t\in T}(X_t-X_{t_0}) \leq O(S) \leq O(1) \sup_{t\in T} \sum_{k\geq 1} 2^{k/2}d(\pi_k(t),\pi_{k-1}(t))
  \]
Since the triangle inequality yields 
  \[
    d(\pi_k(t),\pi_{k-1}(t)) \leq d(t,T_k) + d(t,T_{k-1}) \leq 2d(t,T_{k-1}),
  \]
we obtain the desired generic chaining bound.  
\end{proof}

 \subsection{Concentration of measure \& Applications in the local theory of Banach spaces}
 
In the seventies, the concept of concentration of measure was widely popularized by Vitali Milman as a fundamental phenomenon of study. And indeed, almost fifty years later, this concept has become central both in its theoretical aspects, such as in the geometry of Banach spaces, in random matrix theory, or in the study of random series in Banach spaces, as well as in the world of applied mathematics, for instance, in machine learning, in the theory of compressed sensing, the field of information-based complexity, or even in weather forecasting. The broad scope of concentration of measure techniques and applications is one of the big strengths of this idea. Historically, it originates as a classical consequence of the isoperimetric inequality on the high-dimensional sphere $S^{n-1}:=\{x\in\R^n\,:\,\|x\|_2=1\}$ or on Gauss space, claiming that any Lipschitz function can be seen as a constant function on a set of measure exponentially close to $1$, nowadays referred to as \emph{with overwhelming probability}. 

More precisely, let $f:S^{n-1} \to \R$ be a 1-Lipschitz function and $\varepsilon>0$, then 
  \[
    \nu \left( \theta \in S^{n-1}\,:\, |f(\theta) - M| \ge \varepsilon \right)
     \le 2 \exp(-(n-1) \varepsilon^2),
  \]
where $\nu$ is the uniform measure on the unit sphere and $M$ denotes a median of $f$. In particular, we see that the measure of the set where $f$ is essentially equal to $M$ tends to $1$ exponentially fast as the dimension $n$ of the ambient space tends to infinity. The number $M$ can be replaced by the mean (or expectation) of $f$ over the unit sphere at the price of slight changes in the constants. As famous consequences of the concentration of measure on the sphere one obtains Dvoretzky's theorem on the dimension of almost spherical sections of convex bodies \cite{MR293374, MR445274} and the Johnson--Lindenstrauss lemma \cite{MR0737400} concerning low-distortion embeddings of points from high-dimensional into low-dimensional Euclidean space, which is the fundamental result behind the idea of dimension reduction in data science.

In the Gaussian setting, the concentration bound is in fact a dimension free bound, that is, for any 1-Lipschitz function $F: \R^n \to \R$, 
  \[
    \gamma_n ( |F - \E F(G)| \ge t) \le 2 \exp(-t^2/2),
  \]  
where $\gamma_n$ is the $n$-dimensional Gaussian probability measure and $G$ is a canonical ${\cal N}(0,\mathrm{Id})$ Gaussian random vector in $\R^n$. The development of elegant and sophisticated tools for proving concentration inequalities in various settings has been another major contribution of Talagrand. The publication \cite{mr1361756} is a master piece of his work on this subject. After an extensive presentation of the fundamental concept follows a presentation of various applications, for instance, on bin packing, longest increasing subsequences, longest common subsequences, percolation, chromatic numbers of random graphs, the assignment problem, and sums of vector valued independent random variables. Most of these applications arise in a number of situations such as genetics, speech recognition and, what is very important today, in statistics in high dimensions and machine learning.  Let us present two examples of well-known concentration results of Talagrand. Consider a probability space $(\Omega, \Sigma, \mu)$ and for $N\in\mathbb N$ its $N$-fold product space $(\Omega^N,\mu^N)$; we shall denote by $\Pro$ the corresponding product measure. 
Let $A$ be a subset of $\Omega^N$. For a point $x = (x_1, \ldots, x_N) \in \Omega^N$, the goal is to measure how far $x$ is from the set $A$. A good way to do this is to consider the so-called Hamming distance and, more generally, given a non-negative function $h$ on $\Omega \times \Omega$ such that $h(\omega,\omega) = 0$ for all $\omega \in \Omega$, to define the quantity 
\[
f_h(A,x) := \inf \left\{ \sum_{i \le N} h(x_i, y_i) \,:\, y = (y_1, \ldots, y_N) \in A \right\}.
\]
There is also a different way of measuring how far a point $x$ is from a subset $A$ of $\Omega^N$. Talagrand introduced the set 
\[
U_A(x) := \left \{ (s_i)_{i \le N}\in\{0,1\}^N \,:\, \exists y \in A: (s_i = 0) \implies (x_i = y_i) \right\}.
\]
Roughly speaking, the set $U_A(x)$ encodes to what degree the vector $x$ can be approximated (in the sense of a coordinatewise matching) with a vector from $y\in A$. For instance, if $U_A(x)$ contains many vectors $s$ with small support (i.e., most entries are $0$), then there are many vectors $y$ in $A$ matching $x$ in most coordinates (i.e., on the complement of the support of $s$).
 It is well known that Gustave Choquet, the thesis advisor of Michel Talagrand, gave him maliciously the advice that it is always useful to consider the convex hull. 
 Taking this to heart, let us define $V_A(x)$ to be the convex hull of $U_A(x)$ so that $V_A(x)$ contains $0$ if and only if $x \in A$. Denote by $f_c(A,x)$ the $\ell_2$-distance from zero to $V_A(x)$ and consider for $t\geq 0$ the $t$-enlargement of $A$, which is given by $A_t := \{ x \in \Omega^N\,:\, f_c(A,x) \le t \}$. The following result is commonly referred to as \emph{Talagrand's concentration inequality} and is an isoperimetric-type inequality for product probability spaces representing one of the manifestations of the concentration of measure phenomenon \cite[Theorem 4.1.1]{mr1361756}.
 
\begin{thmalpha}
    For every $A\subseteq \Omega^N$, we have 
    \[
    \int_{\Omega^N} \exp\left( \frac{1}{4} f_c^2 (A,x)\right) d\Pro(x) \le \frac{1}{\Pro(A)}.
    \]
    In particular, 
    \[
    \Pro(A_t) \ge 1 - \frac{1}{\Pro(A)} e^{-t^2/4}.
    \]
\end{thmalpha}
It is not difficult to see that $x \notin A_t$ if and only if 
\[
\forall (\alpha_i)_{i \le N}\in\mathbb R^N \, \exists y \in A:\quad \sum_{i , x_i \ne y_i} \alpha_i \le t \left( \sum_{i=1}^N \alpha_i^2 \right)^{1/2};
\]
see \cite[Lemma 4.1.2]{mr1361756}.
The power of this theorem is illustrated at length in part II of \cite{mr1361756}.  A well known consequence of this result is the following concentration inequality for the supremum of a sum of random vectors \cite[Theorem 13.2]{mr1361756}. 
\begin{thmalpha}
    Consider vectors $(v_i)_{i \le N}$ in a Banach space $W$ and set
    \[
    \sigma := \left( \sup \left\{ \sum_{i \le N} w^* (v_i)^2\,:\, w^* \in W^*, \|w^*\| \le 1 \right\} \right)^{1/2}.
    \]
    Consider a sequence $(Y_i)_{i \le N}$ of independent real valued random variables such that $|Y_i| \le 1$. Denote by $M$ a median of the random variable 
    $\| \sum_{i \le N} Y_i v_i \|$. Then, for every $t>0$, we have
    \[
    \Pro\left( \left| \Big\| \sum_{i \le N} Y_i v_i \Big\| - M \right| \ge t \sigma \right)
    \le 4 \exp\left( - \frac{t^2}{16} \right).
    \]
\end{thmalpha}
A fundamental fact in high dimensions, i.e., when $N$ is large, is that if a subset $A$ is ``large'' then ``most'' of the points of the space are ``close'' to $A$. This is the underlying philosophy of the theory of concentration of measure. ``Large'' and ``most'' are usually understood with respect to a probability measure while there is a huge variety of notions of ``close''. The convex distance introduced by Talagrand is one of them and turned out to be a very appreciated and successful one.

Functional inequalities are at the very heart of various works of Talagrand. The transportation method appeared to be a powerful method to approach the concentration of measure phenomenon, based upon the study of various functional inequalities; some of them are today called ``Talagrand's inequalities''. The transportation cost of a probability measure $\mu$ to a probability measure $\nu$ (both defined on a measurable space $\Omega$) measures how close $\mu$ is to $\nu$ in terms of the ``effort'' required to transport a mass distributed according to $\mu$ into a mass distributed according to $\nu$. To be more precise, let us consider a cost function $w$ on $\Omega \times \Omega$ valued in $\R_+ \cup \{ \infty \}$, symmetric in the arguments, and define the transportation cost $T_w(\mu, \nu)$ as the infimum of the integral
\[
\int_{\Omega \times \Omega} w(x, y) d\pi(x, y)
\]
over all probability measures $\pi$ on $\Omega \times \Omega$ such that $\mu$ is the first marginal of $\pi$ and $\nu$ is the second. \emph{Talagrand's quadratic transportation cost inequality} \cite[Theorem 1.1]{mr1392331} is the following celebrated result.
\begin{thmalpha}
    Consider a measure $\mu$ on $\mathbb R^N$ absolutely continuous with respect to the Gaussian measure $\gamma$ and set $f := \displaystyle \frac{d\mu}{d\gamma}$. Then
    \[
    T_w(\mu, \gamma) \le 2 \int \log f d\mu,
    \]
where the cost function $w$ is given by 
\[
w(x,y) := \|x-y\|_2^2 = \sum_{i \le N} (x_i - y_i)^2.
\]
\end{thmalpha}
The  particular case of the exponential probability measure is also considered in \cite{mr1392331} with a different cost function and leads to sharp concentration inequalities.

To conclude this paragraph, it is important to know that Talagrand's inequalities were mostly motivated by aspects of concentration of measure but also by their applications to the study of the geometry of Banach spaces and probability in Banach spaces. We refer the reader to the wonderful books \cite{LedouxTalagrand1991, Talagrand2014} to get an overview of the variety and the numerous applications in Banach space theory. A paper that we particularly like is \cite{mr1368249} where Talagrand studied in great detail the sections of smooth convex bodies. When you do not know how to prove that a subset of your vectors (or functions) satisfies a particular property, it is particularly useful to select them at random and to prove that with positive (very often ``high'') probability, the choice of the random element will satisfy your desired property. Bourgain used in a huge variety of problems in harmonic analysis the so called  ``method of selector''. In his solution to the $\Lambda(p)$ problem, Bourgain developed a masterpiece of a technique and used different arguments for $2<p<3$, $3<p<4$, and $p>4$. The study of particular majorizing measures allowed Talagrand to establish that, in fact, the key point to solve the $\Lambda(p)$ problem is that $L_p$ spaces are 2-smooth for $p>2$, which means that for some constant $C\in(0,\infty)$,  for every $f \in L_p$ with $\|f\| = 1$ and every $g \in L_p$ with $\|g\| \le 1$, we have
\[
\frac{1}{2} \left(\left \|f+g \right\|  +  \left \| f-g \right\| \right)
\le
1  + C \| g\|^2.
\]
This inequality holds true in $L_p$ spaces for $p \ge 2$ and the constant $C$ can be chosen to be $(p-1)/2$.
There are various papers of Talagrand about selection of characters that are based on a deep construction of a majorizing measure. Roughly speaking, in a 2-smooth Banach space $X$, for a family of vectors $x_1, \ldots, x_N$, we associate the set ${\cal F} \subseteq \R^N$ 
\[
{\cal F} := \left\{  x^* (x_i)^2 \,:\, x^* \in X^*, \|x^*\| \le 1  \right\}.
\]
The goal is to select a subset $I$ of $\{1, \ldots, N\}$ such that the family $(x_i)_{i \in I}$ satisfies a particular property. To do that, you define a family of independent identically distributed random variables $\delta_1, \ldots, \delta_N$ taking values  in $\{0, 1\}$ with mean $\delta$.  The goal is to bound the following type of process
\[
\E \sup_{f \in {\cal F}} \sum_{i=1}^N \delta_i f_i
\]
choosing ``correctly'' the parameter $\delta$, and this is where majorizing measures enter the game, see \cite[Theorem 1.3]{mr1368249}.

\section{Talagrand's Work -- Spin glasses and the Parisi formula}
 
 In the field of physics dealing with the macroscopic and microscopic physical properties of matter, a \emph{spin glass} is a magnetic state characterized by randomness (besides cooperative behavior in freezing of spins at a temperature called ``freezing temperature'') and is therefore also referred to as a ``disordered'' magnetic state in contrast to a ferromagnet where magnetic spins are ordered, i.e., all align in the same direction. The term ``glass'' is meant to draw an analogy between the magnetic disorder and the highly irregular atomic bond structure in amorphous materials such as window glass.  
 
To be able to get an idea of Talagrand's contribution in the theory of spin glasses, let us start with the outstanding work of theoretical physicist Giorgio Parisi (born 4 August 1948), whose research focuses on quantum field theory, statistical mechanics, and complex systems. One particular and famous result is his exact (but not rigorous) solution of the Sherrington--Kirkpatrick (SK) model of spin glasses \cite{Parisi1980}, introduced by David Sherrington and Scott Kirkpatrick in 1975 \cite{SherringtonKirkpatrick1975} believing that their model was amenable to more straightforward mathematical analysis than the more realistic lattice model of Samuel F. Edwards and Phil W. Anderson \cite{EdwardsAnderson1975} (a simple extension of the well known nearest neighbor Ising model). The SK model is nothing else than an Ising model, i.e., a mathematical model of ferromagnetism in statistical mechanics, with long range frustrated ferro- as well as antiferromagnetic bonds and corresponds to a mean-field approximation (i.e., the reference system minimizing the right-hand side in Bogoliubov's inequality estimating the free energy of the original Hamiltonian of the model) of spin glasses describing the slow dynamics of the magnetization and the complex non-ergodic equilibrium state.
The SK model, even though a purely mathematical object, has been studied by Parisi and followers by methods not recognized as legitimate by most mathematicians, i.e., arguments were rather heuristic and built on physical intuition; nonetheless, it is important to note that the Parisi solution is presented in form of an ingenious and clever Ansatz. In his groundbreaking paper \cite{Talagrand2006}, Talagrand corrected this discrepancy and made one of the central predictions of Parisi (which were instrumental in awarding the 2021 Nobel Prize in Physics to him) rigorous, the computation of the ``free energy'' of the SK model. Talagrand's proof uses the so-called Guerra's interpolation scheme, a very powerful method based on a comparison and interpolation argument on sets of Gaussian random variables (we refer to \cite{Guerra2005} for details).

Let us recall that the random Hamiltonian of the SK model for spin glasses is given by the mean field expression (see, e.g., \cite{SherringtonKirkpatrick1975})
  \[
    H_N(\sigma) := -\frac{1}{\sqrt{N}} \sum_{1\leq i<j \leq N} g_{ij}\sigma_i\sigma_j - h\sum_{i=1}^N\sigma_i,
  \]
where $\sigma=(\sigma_1,\dots,\sigma_N) \in \Sigma_N:=\{-1,1\}^N$ are the spin configurations and $(g_{ij})_{i<j}$ are independent and identically distributed standard Gaussians, the so-called random couplings or external quenched disorder. The first term on the right-hand side of the Hamiltonian represents a long range random two body interaction, while the second is an external component representing the interaction of the spins with a fixed external magnetic field $h$. The disorder dependent normalization constant at inverse temperature $\beta\geq 0$, known as partition function, is defined as
  \[
    Z_N := Z_N(\beta):= \sum_{\sigma\in \Sigma_N} e^{- \beta H_N(\sigma)}
  \]
and allows to introduce the corresponding configuration probability given by the random Gibbs measure
  \[
    \mu_{\beta}(\sigma) = \frac{e^{-\beta H_N(\sigma)}}{Z_N}
  \]
on the space of spin configurations.  
The quantity of interest and at the heart of Talagrand's contribution is (essentially)
  \[
    f_N := -\frac{1}{\beta N}\E\big[\log Z_N \big]
  \]
and called (the quenched average of) the free energy per site; expectation is always understood in the disorder, i.e., in the randomness of the Hamiltonian. Let us note that $\beta \mapsto \frac{1}{N} \log Z_N $ is a (random) convex function and, consequently, $\beta \mapsto \frac{1}{N}\E\big[\log Z_N \big]$ a convex function.  

In a complex proof, Talagrand managed to compute the asymptotic free energy of the SK model (actually a function closely related to it is computed). A rigorous lower bound had been obtained before by Guerra \cite{Guerra2002} by means of an interpolation scheme, a technique that is also the backbone of Talagrand's proof, while the mere existence of the (infinite volume) limits had been proven by Guerra and Toninelli in \cite{GuerraToninelli2002}. 

The groundbreaking result of Talagrand is the following theorem (see \cite[Theorem 1.1]{Talagrand2006}). 

\begin{thmalpha}[The Parisi formula]\label{thm:parisi formula}
We have that
  \[
    \lim_{N\to\infty} \frac{1}{N} \E \big[\log Z_N \big] = P(\xi,h),
  \]   
where the quantity $P$ depends on the function $\xi(x):=\beta x^2/2$ and the strength $h$ of the external magnetic field. 
\end{thmalpha}

The Parisi conjecture, which is confirmed by Theorem \ref{thm:parisi formula}, was arguably the most famous open problem in the theory of spin glasses, but there are still quite a number of fundamental predictions of the Parisi theory that remain open problems to the present day and we refer the reader to, e.g., \cite{Talagrand2007} for a first overview. We shall not immerse in the delicate and technical proof or even try to present its key elements, and rather refer the interested reader to the original paper \cite{Talagrand2006} and the references cited therein. In particular, in the introduction the author nicely elaborates on the developments that lead to his proof.

\subsection*{Acknowledgement}
Joscha Prochno's research is supported by the German Research Foundation (DFG) under project 516672205.

\bibliographystyle{plain}
\bibliography{abel_talagrand, annee-biblio-talagrand}

\bigskip
\bigskip

	\medskip
	
	\small
	
	\noindent \textsc{Olivier Gu\'edon:} Laboratoire d'Analyse et de Math\'ematiques Appliqu\'ees, Univ Gustave Eiffel, Univ Paris Est Creteil, CNRS, LAMA UMR8050 F-77447 Marne-la-Vallée, France.
	
	\noindent
	{\it E-mail:} \texttt{olivier.guedon@univ-eiffel.fr}
	
	\medskip

	\noindent \textsc{Joscha Prochno:} Faculty of Computer Science and Mathematics, University of Passau, Dr.-Hans-Kapfinger-Str. 30, 94032 Passau, Germany.
	
	\noindent
	{\it E-mail:} \texttt{joscha.prochno@uni-passau.de}

\end{document}